\documentclass[]{interact}

\usepackage{epstopdf}
\usepackage[caption=false]{subfig}
\usepackage[shortlabels]{enumitem}

\usepackage[numbers,sort&compress]{natbib}
\bibpunct[, ]{[}{]}{,}{n}{,}{,}
\makeatletter
\def\NAT@def@citea{\def\@citea{\NAT@separator}}
\makeatother

\theoremstyle{plain}
\newtheorem{theorem}{Theorem}[section]

\newtheorem{conjecture}{Conjecture}

\theoremstyle{definition}
\newtheorem{definition}[theorem]{Definition}

\theoremstyle{remark}

\newtheorem{notation}{Notation}

\begin{document}


\title{The permanent and diagonal products on the set of nonnegative matrices with bounded rank}

\author{
\name{Yair Lavi\thanks{CONTACT Y. Lavi. Email: yairlv@gmail.com}}
}

\maketitle

\begin{abstract}
We formulate conjectures regarding the maximum value and maximizing matrices of the permanent and of diagonal products on the set of stochastic matrices with bounded rank. We formulate equivalent conjectures on upper bounds for these functions for nonnegative matrices based on their rank, row sums and column sums. 
In particular we conjecture that the permanent of a singular nonnegative matrix is bounded by $\frac{1}{2}$ times the minimum of the product of its row sums and the product of its column sums, and that the product of the elements of any diagonal of a singular nonnegative matrix is bounded by $\frac{1}{4}$ times the minimum of the product of its row sums and the product of its column sums.
\end{abstract}


\section{Introduction}

The permanent of a real $n\times n$ matrix is defined by $per(A)=\sum_{\sigma\in S_{n}}\prod_{i=1}^{n} a_{i\sigma(i)}$
 
Consider the familiar inequality for the permanent of an $n \times n$ (row or column) stochastic matrix $A$ \cite{per_minc}: 

\begin{equation}
per(A)\leq 1
\label{eq:st_inequality}
\end{equation}

 Equality is obtained iff $A$ is a permutation matrix, and in particular a necessary condition for equality is that $rank(A)=n$.
 One could naturally ask then: \\What is the maximum value of the permanent function and what are the maximizing matrices in the set of stochastic matrices with rank bounded from above by some $k<n$ ?\\
In this note we conjecture an answer to this question.

An equivalent result to inequality \ref{eq:st_inequality} is the following bound for the permanent of a nonnegative matrix $A$:

\begin{equation}
per(A)\leq min(\prod_{i=1}^{n}r_{i},\prod_{i=1}^{n}c_{i})
\label{eq:nn_inequality}
\end{equation}

where $(r_{1},\ldots ,r_{n})$ and $(c_{1},\ldots ,c_{n})$ are the row sums and column sums of $A$  respectively.

We formulate our conjecture in the equivalent "nonnegative formulation"  as well.\par

Let $\sigma\in S_{n}$ be any permutation. 
We call $\prod_{i=1}^{n} a_{i\sigma(i)}$ a diagonal product of $A$.\\
Similarly to the permanent conjecture we formulate a conjecture for the maximum of any diagonal product on the set of stochastic matrices with bounded rank and also present it in the equivalent "nonnegative formulation".

We begin in Section 2 with notations and definitions. In Section 3, we formulate the conjecture regarding the permanent function and in Section 4 we formulate the conjecture regarding the diagonal products.

\section{Preliminaries}

We denote by $\mathcal{M}_{n}$ the set of real $n\times n$  matrices
and by $\mathcal{M}_{n}^{+}$ the set of nonnegative matrices in $\mathcal{M}_{n}$. \\

We denote the set of nonnegative matrices of rank less than or equal to k by $$M_{n}^{+\overline{k}}:=\{A\in \mathcal{M}_{n}^{+}|  rank(A)\leq k\}$$

Let $R_{n}$ denote the set of row stochastic matrices in $\mathcal{M}_{n}^{+}$, so $R_{n}$ is the set of $n\times n$  real matrices A =($a_{ij}$) satisfying:
$$\sum_{j=1}^{n}a_{ij}=1 \quad  i=1,\ldots ,n$$ and  $$0\leq a_{ij}\quad  1\leq i,j\leq n   $$we denote the set of row stochastic matrices of rank less than or equal to k  by $$R_{n}^{\overline{k}}:=\{A\in R_{n}|  rank(A)\leq k\}$$
\\
 
 Similarly, we will write  $L_{n}$, and $L_{n}^{\overline{k}}$ for the sets of column stochastic matrices and column stochastic matrices of rank less than or equal to k, respectively.

For a matrix $A\in \mathcal{M}_{n}^{+}$ with rows $A_{1},\ldots,A_{n}$ and row sums $r_{i}>0,\quad  i=1,\ldots ,n$, we will denote by $\overline{A}^{r}$ the row stochastic matrix whose rows are $\frac{A_{1}}{r_{1}},\frac{A_{2}}{r_{2}},\ldots ,\frac{A_{n}}{r_{n}}$.
Similarly, for a matrix $A\in \mathcal{M}_{n}^{+}$ with columns $\tilde{A}_{1},\ldots ,\tilde{A}_{n}$ and column sums $c_{j}>0, \quad j=1,\ldots, n$, we will denote by $\overline{A}^{c}$ the column stochastic matrix whose columns are $\frac{\tilde{A}_{1}}{c_{1}},\frac{\tilde{A}_{2}}{c_{2}},\ldots ,\frac{\tilde{A}_{n}}{c_{n}}$.

\begin{notation}
	[\textbf{The Set of Compositions of n with s parts}]
	We denote by $P_{n;s}\subset\mathbb{R}^{s}$ the set of all positive integer compositions of n, with exactly s parts, so $\vec{r}\in P_{n;s}$ iff $\sum_{i=1}^{s} r_{i}=n$ ,  $r_{i}\in\mathbb{N}, i=1,2\ldots ,s$.
\end{notation}

Let $J_{r}$ be the $r\times r$ doubly stochastic matrix whose entries all equal to $\frac{1}{r}$.
\begin{definition}
	[\textbf{ s-composition matrix}]
	Given a composition $\vec{r}=(r_{1},r_{2},..,r_{s})\in P_{n;s}$ we shall call the doubly stochastic matrix $J_{r_{1}}\oplus J_{r_{2}}.....\oplus J_{r_{s}} $ an \textit{s-composition matrix} and denote it by $J_{\vec{r}}$ 
	
\end{definition}
 
 \begin{notation}
 	Let $\pi:S_{n}\to\mathcal{M}_{n}$ be the natural permutation representation of the symmetric group $S_{n}$ defined by 
 	$\pi(\sigma)_{ij}=\delta_{\sigma(i),j}$ and denote its image, the set of permutation matrices, by $\mathfrak{S}_{n}=\pi(S_{n})$. 
 \end{notation}

\section{Conjectures for the permanent function}

We formulate our conjecture in two equivalent forms, starting with the stochastic matrices formulation:

\begin{conjecture}  \label{conj:per}
	\textbf{Maximum of the permanent function on the set of stochastic matrices with bounded rank.} For $k,n\in\mathbb{N}$, $1\leq k\leq n$, let $R_{n}^{\overline{k}}$ ($L_{n}^{\overline{k}}$) be the set of $n\times n$ row (column) - stochastic matrices of rank less than or equal to $k$, and let $r,s\in\mathbb{N}$ be the unique integers such that $n=rk+s$, $0\leq s<k$  .Then  $$A\in R_{n}^{\overline{k}}\cup L_{n}^{\overline{k}}\implies per(A)\leq (\frac{r!}{r^{r}})^{k-s}\times(\frac{(r+1)!}{(r+1)^{r+1}})^{s}$$  with equality iff  $A=PJ_{\vec{r}}Q$ where $\vec{r}=(\underbrace{r,r,\ldots,r}_{(k-s)\text{-times}},\underbrace{r+1,r+1,\ldots,r+1}_{s\text{-times}})$ and $P,Q\in\mathfrak{S}_{n}$ are any two permutation matrices. 
	
	In particular the maximizing matrices are doubly stochastic matrices of $rank=k$.
	
\end{conjecture}

We now present the equivalent formulation for nonnegative matrices.
The equivalence of the two formulations is trivial.

\begin{conjecture} \label{conj:per_nn}
	\textbf{Permanent function of nonnegative matrices with bounded rank.}
	For $k,n\in\mathbb{N}$, $1\leq k\leq n$, let $M_{n}^{+\overline{k}}$ be the set of $n\times n$ nonnegative matrices of rank less than or equal to $k$, let $r,s\in\mathbb{N}$ be the unique integers such that $n=rk+s$, $0\leq s<k$ and let  $\vec{r}=(\underbrace{r,r,\ldots,r}_{(k-s)\text{-times}},\underbrace{r+1,r+1,\ldots,r+1}_{s\text{-times}})$ .Then  $$A\in M_{n}^{+\overline{k}}\implies per(A)\leq min(\prod_{i=1}^{n}r_{i},\prod_{i=1}^{n}c_{i})(\frac{r!}{r^{r}})^{k-s}\times(\frac{(r+1)!}{(r+1)^{r+1}})^{s}$$ where $(r_{1},\ldots ,r_{n})$ and $(c_{1},\ldots ,c_{n})$ are the row sums and column sums of $A$ respectively, 
	
	with equality iff one of the following holds:
	\begin{enumerate}
		\item  $A$ has a zero row or a zero column 
		\item  $0<\prod_{i=1}^{n}r_{i}< \prod_{i=1}^{n}c_{i}$ and  $\overline{A}^{r}=PJ_{\vec{r}}Q$ for some $P,Q\in\mathfrak{S}_{n}$ 
		\item  $0<\prod_{i=1}^{n}c_{i}< \prod_{i=1}^{n}r_{i}$ and  $\overline{A}^{c}=PJ_{\vec{r}}Q$ for some $P,Q\in\mathfrak{S}_{n}$ 
		\item  $0< \prod_{i=1}^{n}r_{i}=\prod_{i=1}^{n}c_{i}$ and  $\overline{A}^{r}=PJ_{\vec{r}}Q$ and $\overline{A}^{c}=P'J_{\vec{r}}Q'$ for some $P,Q,P',Q'\in\mathfrak{S}_{n}$
	\end{enumerate}

\end{conjecture}
\pagebreak
In particular for singular nonnegative matrices we make the following conjecture:
\begin{conjecture}
\textbf{Permanent function of nonnegative singular matrices.}
Let $A$ be a singular nonnegative matrix and let  $\vec{r}=(\underbrace{1,1,\ldots,1}_{(n-2)\text{-times}},2)$.
then $$ per(A)\leq \frac{min(\prod_{i=1}^{n}r_{i},\prod_{i=1}^{n}c_{i})}{2}$$ where $(r_{1},\ldots ,r_{n})$ and $(c_{1},\ldots ,c_{n})$ are the row sums and column sums of $A$ respectively,

 with equality iff one of the following holds:
 \begin{enumerate}
 	\item $A$ has a zero row or a zero column 
 	\item  $0<\prod_{i=1}^{n}r_{i}< \prod_{i=1}^{n}c_{i}$ and  $\overline{A}^{r}=PJ_{\vec{r}}Q$ for some $P,Q\in\mathfrak{S}_{n}$ 
 	\item  $0<\prod_{i=1}^{n}c_{i}< \prod_{i=1}^{n}r_{i}$ and  $\overline{A}^{c}=PJ_{\vec{r}}Q$ for some $P,Q\in\mathfrak{S}_{n}$ 
 	\item  $0< \prod_{i=1}^{n}r_{i}=\prod_{i=1}^{n}c_{i}$ and  $\overline{A}^{r}=PJ_{\vec{r}}Q$ and $\overline{A}^{c}=P'J_{\vec{r}}Q'$ for some $P,Q,P',Q'\in\mathfrak{S}_{n}$
 \end{enumerate}

\end{conjecture}

\section{Conjectures for diagonal products}

We formulate our conjecture in two equivalent forms, starting with the stochastic matrices formulation:

\begin{conjecture} \label{conj:diag}
	\textbf{Maximum of diagonal products on the set of stochastic matrices with bounded rank.} 
	Let $\sigma\in S_{n}$ be any permutation. For $k,n\in\mathbb{N}$, $1\leq k\leq n$, let $R_{n}^{\overline{k}}$ ($L_{n}^{\overline{k}}$) be the set of $n\times n$ row (column) - stochastic matrices of rank less or equal to $k$, and let $r,s\in\mathbb{N}$ be the unique integers such that $n=rk+s$, $0\leq s<k$  .Then  $$A\in R_{n}^{\overline{k}}\cup L_{n}^{\overline{k}}\implies \prod_{i=1}^{n}a_{i\sigma(i)}\leq (\frac{1}{r^{r}})^{k-s}\times(\frac{1}{(r+1)^{r+1}})^{s}$$ with equality iff  $A=P^{t}J_{\vec{r}}P\pi (\sigma)$ where $\vec{r}=(\underbrace{r,r,\ldots,r}_{(k-s)\text{-times}},\underbrace{r+1,r+1,\ldots,r+1}_{s\text{-times}})$ and $P\in\mathfrak{S}_{n}$ is any permutation matrix.

\end{conjecture}

We now present the equivalent formulation for nonnegative matrices.
The equivalence of the two formulations is trivial.

\begin{conjecture} \label{conj:any_diag_nn}
	\textbf{Diagonal products of nonnegative matrices with bounded rank.} 
	
	Let $\sigma\in S_{n}$ be any permutation. 
	For $k,n\in\mathbb{N}$, $1\leq k\leq n$, let $M_{n}^{+\overline{k}}$ be the set of $n\times n$ nonnegative matrices of rank less or equal to $k$, let $r,s\in\mathbb{N}$ be the unique integers such that $n=rk+s$, $0\leq s<k$ and let $\vec{r}=(\underbrace{r,r,\ldots,r}_{(k-s)\text{-times}},\underbrace{r+1,r+1,\ldots,r+1}_{s\text{-times}})$  Then  $$A\in M_{n}^{+\overline{k}}\implies \prod_{i=1}^{n}a_{i\sigma(i)}\leq min(\prod_{i=1}^{n}r_{i},\prod_{i=1}^{n}c_{i}) (\frac{1}{r^{r}})^{k-s}\times(\frac{1}{(r+1)^{r+1}})^{s}$$ 
	\\
	where $(r_{1},\ldots ,r_{n})$ and $(c_{1},\ldots ,c_{n})$ are the row sums and column sums of $A$  respectively,
	with equality iff one of the following holds:
	\begin{enumerate}
		\item $A$ has a zero row or a zero column 
		\item  $0<\prod_{i=1}^{n}r_{i}< \prod_{i=1}^{n}c_{i}$ and  $\overline{A}^{r}=P^{t}J_{\vec{r}}P\pi(\sigma)$ for some $P\in\mathfrak{S}_{n}$ 
		\item  $0<\prod_{i=1}^{n}c_{i}< \prod_{i=1}^{n}r_{i}$ and  $\overline{A}^{c}=P^{t}J_{\vec{r}}P\pi(\sigma)$ for some $P\in\mathfrak{S}_{n}$ 
		\item  $0< \prod_{i=1}^{n}r_{i}=\prod_{i=1}^{n}c_{i}$ and  $\overline{A}^{r}=P^{t}J_{\vec{r}}P\pi(\sigma)$ and $\overline{A}^{c}=P'^{t}J_{\vec{r}}P'\pi(\sigma)$ for some $P,P'\in\mathfrak{S}_{n}$
	\end{enumerate}

\end{conjecture}

In particular for singular nonnegative matrices we make the following conjecture:
\begin{conjecture}
	\textbf{Diagonal products of nonnegative singular matrices.}
	Let $\sigma\in S_{n}$ be any permutation, let $A$ be a singular nonnegative matrix and let  $\vec{r}=(\underbrace{1,1,\ldots,1}_{(n-2)\text{-times}},2)$
	then $$ \prod_{i=1}^{n}a_{i\sigma(i)}\leq \frac{min(\prod_{i=1}^{n}r_{i},\prod_{i=1}^{n}c_{i})}{4}$$ where $(r_{1},\ldots ,r_{n})$ and $(c_{1},\ldots ,c_{n})$ are the row sums and column sums of $A$ respectively, 
	with equality iff one of the following holds:
	\begin{enumerate}
	\item $A$ has a zero row or a zero column 
	\item  $0<\prod_{i=1}^{n}r_{i}< \prod_{i=1}^{n}c_{i}$ and  $\overline{A}^{r}=P^{t}J_{\vec{r}}P\pi(\sigma)$ for some $P\in\mathfrak{S}_{n}$ 
	\item  $0<\prod_{i=1}^{n}c_{i}< \prod_{i=1}^{n}r_{i}$ and  $\overline{A}^{c}=P^{t}J_{\vec{r}}P\pi(\sigma)$ for some $P\in\mathfrak{S}_{n}$ 
	\item  $0< \prod_{i=1}^{n}r_{i}=\prod_{i=1}^{n}c_{i}$ and  $\overline{A}^{r}=P^{t}J_{\vec{r}}P\pi(\sigma)$ and $\overline{A}^{c}=P'^{t}J_{\vec{r}}P'\pi(\sigma)$ for some $P,P'\in\mathfrak{S}_{n}$
\end{enumerate}

\end{conjecture}



\appendix

\end{document}